\documentclass[11pt,onecolumn]{article}

\usepackage{graphicx}      % include this line if your document contains figures
\usepackage{algorithm, algorithmic, setspace}
\usepackage{graphicx} % for pdf, bitmapped graphics files
\usepackage{amsmath} % assumes amsmath package installed
\usepackage{amssymb}  % assumes amsmath package installed
\usepackage{amsthm}
\usepackage{amsfonts}
\usepackage{color}
\usepackage[english]{babel}

\newcommand{\argmin}[1]{\underset{#1}{\mathrm{argmin\,}}}

\newcommand{\xtrue}{\widetilde{x}}
\newcommand{\atrue}{\widetilde{a}}

\newcommand{\sign}{\text{sign}}
\newcommand{\xstar}{x^{\star}}
\newcommand{\astar}{a^{\star}}
\newcommand{\R}{\mathbb{R}}

\newcommand{\soft}{\mathrm{S}}

\newcommand{\su}{\mathcal{S}}
\newcommand{\suc}{\bar{\mathcal{S}}}

\newcommand{\Obs}{\mathcal{O}}

\newcommand{\Is}{I_{\su}}
\newcommand{\Isc}{I_{\suc}}
\newcommand{\irr}{\rho}

\newcommand{\rank}{\mathrm{rank}}

\newcommand{\xhat}{\hat{x}}
\newcommand{\ahat}{\hat{\mathbf{a}}}
\newcommand{\yhat}{\hat{\mathbf{y}}}

\newtheorem{problem}{Problem}
\newtheorem{theorem}{Theorem}

\newtheorem{example}{Example}

\begin{document}
\title{Lasso-based state estimation for\\ cyber-physical systems under sensor attacks}
\author{V. Cerone\thanks{The authors are with the Department of Control and Computer Engineering, Politecnico di Torino, Italy (e-mail: sophie.fosson@polito.it). This paper is part of the project NODES which has received funding from the MUR – M4C2 1.5 of PNRR with grant agreement no. ECS00000036.}, $~~~$   S. M. Fosson\footnotemark[1],  $~~~$  D. Regruto\footnotemark[1],  $~~~$  F. Ripa\footnotemark[1]}

\maketitle

%
%\address[Poli]{Department of Control and Computer Engineering,\\ Politecnico di Torino, Italy (e-mail: sophie.fosson@polito.it).}
%
\begin{abstract}The development of algorithms for secure state estimation in vulnerable cyber-physical systems has been gaining attention in the last years. A consolidated assumption is that an adversary can tamper a relatively small number of sensors. In the literature, block-sparsity methods exploit this prior information to recover the attack locations and the state of the system.
In this paper, we propose an alternative, Lasso-based approach and we analyse its effectiveness. In particular, we theoretically derive conditions that guarantee successful attack/state recovery, independently of established time sparsity patterns. Furthermore, we develop a sparse state observer, by starting from the iterative soft thresholding algorithm for Lasso, to perform online estimation.
Through several numerical experiments, we compare the proposed methods to the state-of-the-art algorithms.
\end{abstract}
%
%\begin{keyword}
%Cyber-physical systems, sensor attacks, secure state estimation, sparse optimization, Lasso, Luenberger-like observers
%\end{keyword}

%\end{frontmatter}
%===============================================================================

\section{Introduction}\label{sec:IN}
% \begin{enumerate}
%  \item Cyber
% \end{enumerate}
%
A cyber-physical system (CPS) is a collection of computing devices that interact with the physical world, through sensors and actuators, and with one another, through communication networks. Applications of the CPS paradigm include industrial control processes, smart power grids, wireless sensor networks, electric ground vehicles and cooperative driving technologies.

The distributed nature of CPSs is ambivalent in terms of security: on the one hand, it is more resilient to faults with respect to centralized systems; on the other hand, it is exposed to adversaries, either in terms of physical access to sensors or cyber access to data transmission networks. Examples range from non-invasive spoofing physical attacks,  as illustrated by \cite{sho13}, to cyber attacks to SCADA systems as done by the Stuxnet worm; see \cite{lan11}. %We refer the reader to the survey \cite{duo22} for an extended overview on attacks in CPSs.
Since CPSs act on the physical world, cyber-physical attacks layer may yield serious consequences to physical processes and to human beings, which makes the security problem very critical.

In the last decade, a substantial work has focused on the security of CPSs. A relevant research line considers the problem of secure state estimation (SSE) for CPSs in the presence of sensor attacks, that inject false data to manipulate the measurements.
As a matter of fact, a malicious injection perturbs the data as a noise. However, we expect that an adversary conceives an unpredictable intrusion, that is, we have no information on its dynamics, boundedness and probabilistic description. The unique realistic assumption on sensor attacks is sparsity: only a relatively small number of sensors is accessible, due to, e.g., large dimensionality, physical deployment and sensor heterogeneity of CPSs.

In this paper, we consider CPSs described by discrete-time (DT) linear time-invariant (LTI) dynamical systems, with sparse sensor attacks. The identification of the attack support, i.e., the subset of tampered sensors, is a combinatorial problem, which does not scale well for large dimensional systems. However, by leveraging the sparsity assumption, one can exploit $\ell_1$-based sparsity-promoting decoders to recast the problem into constrained convex optimization; see, e.g.,  \cite{faw14,paj17}. Since these approaches are still computationally intense, \cite{sho16} introduce a faster event-triggered projected gradient (ETPG) approach, whose structure is prone to recursive SSE. The provided sufficient conditions for the convergence of ETPG are quite restrictive. \cite{sho17} address this issue by a satisfiability modulo theory approach, called Imhotep-SMT, which returns the exact solution in short time for problems of small/medium dimensions. However, Imhotep-SMT is combinatorial, thus critical for large-scale problems.

%The approaches proposed by \cite{faw14,sho16,sho17} mainly focus on noise-free problems; they may envisage the presence of noise, beyond the attacks, but with few robustness guarantees. %Moreover, these works assume either the time invariance of the attack support or the knowledge of a tight upper bound of the number of attacks.

In this paper, we propose a different approach to SSE of CPSs under sparse sensor attacks, based on Lasso; see \cite{tib96}. Specifically, we define a Lasso formulation and we analyse its effectiveness.

The proposed Lasso-based model builds on $\ell_1$ relaxation, but, differently from \cite{faw14,paj17}, it gives rise to an unconstrained optimization problem, which can be solved through low-complex recursive algorithms. By elaborating on this point, the second contribution of the paper is a recursive SSE method exploiting new data as soon as they become available, with the final aim of performing online SSE. More precisely, we design a sparsity-promoting Luenberger-like observer by starting from the iterative soft thresholding algorithm for Lasso.

Finally, we propose numerical experiments that show the effectiveness of the proposed methods with respect to the state-of-the-art approaches, in terms of estimation accuracy and execution time.

%to check its batch state estimation performance. While in \cite{faw14,sho16,sho17} the focus is on the exact state estimation in presence of attacks, our point is on the correct detection of the attack support; indeed,if we know the support, according to standard observability assumption, one can recover the state from the attack-free measurements. To understand the effectiveness of Lasso to recover the attack support, we develop an irrepresentable condition analysis specialized for CPSs under sensors attacks, which provides interpretable conditions for success.
%
%As in \cite{faw14,sho16}, in this paper,  we focus on noise-free models, although the proposed approach is not strictly tied to these features.

We organize the paper as follows. In Sec. \ref{sec:PS}, we state the problem and we illustrate the  background. In Sec. \ref{sec:PA}, we introduce the proposed Lasso approach and we theoretically analyse it. In Sec. \ref{sec:OB}, we extend the method to recursive and online SSE, by developing a state observer. Finally, we devote Sec. \ref{sec:NR} to numerical experiments and we draw some conclusions in Sec. \ref{sec:C}.
\section{Problem Statement}\label{sec:PS}
By following \cite{faw14,sho16}, we consider CPSs that can be modeled as DT LTI dynamical systems
\begin{equation}\label{system}
 \begin{split}
  x(k+1)&=Ax(k)\\
  y(k)&=Cx(k)+a(k)
 \end{split}
\end{equation}
where  $x(k)\in\R^n$ is the state, $y(k)\in\R^p$ is the measurement vector, $a(k)\in\R^p$ is the attack vector, $A\in\R^{n,n}$ and $C\in\R^{p,n}$. We assume that each sensor $i$ takes a measurement $y_i(k)$; if $a_i(k)\neq 0$, sensor $i$ is under attack. We assume that $a(k)$ is sparse, i.e., few sensors are under attack at each $k$. %the attack support is time-invariant and has cardinality $s\ll p$.
Since the presence of a known input  does not impact on the formulation of the problem, see \cite{sho16}, for simplicity of notation we consider a zero-input model. Moreover, as in \cite{faw14,sho16}, we consider a noise-free model for our theoretical analysis, while we envisage measurement noise in some numerical experiments.

The SSE problem is as follows.
\begin{problem}\label{p1}
For some $\tau\leq n$ and $k\geq \tau-1$ given $A$, $C$ and $y=(y(k-\tau+1)^\top, \dots,  y(k)^\top)^\top\in\R^{p\tau}$, estimate the $\tau$-delayed state $x(k-\tau+1)$  in the presence of sparse sensor attacks.
\end{problem}
% For some $\tau\leq n$, given $A$, $C$ and $y=(y(0)^\top, \dots,  y(\tau-1)^\top)^\top\in\R^{p\tau}$, we aim at estimating $x(0)$ in the presence of sparse sensor attacks $a(k)$.
% \end{problem}
% \begin{remark}\label{rem}
% %Problem \ref{p1} is equivalent to say that, at any time instant $k$, we aim at estimating the $\tau$-delayed state $x(k-\tau+1)$ from $(y(k-\tau+1)^\top, \dots,  y(k)^\top)^\top$.
% \end{remark}

%To estimate $x(0)$, we need to estimate $a(0)$,  $a(1)$, $\dots$, $a(\tau-1)$.
Let us denote $\atrue=(a(k-\tau+1)^\top, \dots,  a(k)^\top)^\top\in\R^{p\tau}$ and $\xtrue=x(k-\tau+1)\in\R^n$, while $I\in\{0,1\}^{p\tau,p\tau}$ is the identity matrix.
From \eqref{system}, we have
\begin{equation}\label{fact}
y =\begin{pmatrix}\Obs & I
 \end{pmatrix}\begin{pmatrix}
               \xtrue\\
               \atrue
              \end{pmatrix}, \text{ where } \Obs=\begin{pmatrix}
C\\CA \\ \vdots \\ CA^{\tau-1}
\end{pmatrix}\in\R^{p\tau,n}
\end{equation}
%where
%$y=\begin{pmatrix}
%  y(0)\\
%  \vdots\\
%  y(\tau-1)
% \end{pmatrix}\in\R^{p\tau}$
% $ \Obs=\begin{pmatrix}
% C\\CA \\ \vdots \\ CA^{\tau-1}
% \end{pmatrix}\in\R^{p\tau,n}$,
%  $\atrue=\begin{pmatrix}
%                a(k-\tau+1)\\
%                \vdots\\
%                a(k)
%               \end{pmatrix}\in\R^{p\tau}$,

%
If $\tau=n$, $\Obs$ is the observability matrix of the attack-free system. We assume that the attack-free system is observable, i.e., $\rank(\Obs)=n$.
In principle, we can estimate $(\xtrue, \atrue)$ by solving
\begin{equation}\label{cs}
 y=\begin{pmatrix}\Obs & I
 \end{pmatrix}\begin{pmatrix}
               x\\
               a
              \end{pmatrix}
\end{equation}
in the variables $x\in\R^n$, $a\in\R^{p\tau}$. Nevertheless, since $\begin{pmatrix}\Obs & I
 \end{pmatrix}\in\R^{p\tau,n+p\tau}$, the system is inherently underdetermined. However, by taking into account the sparsity of $\atrue$, we can exploit the compressed sensing theory to find a sparse solution to \eqref{cs}; see, e.g., \cite{fou13}. Some works assume that the attack support is time-invariant, with cardinality $s\ll p$, i.e., $\atrue$ is block-sparse. In particular, \cite{faw14} exploit an $\ell_1$/$\ell_r$ norm approach for block-sparse signal recovery, while \cite{sho16} develop a block-based hard thresholding algorithm, that alternates gradient descent and event-triggered projection onto $\R^n\times\mathbb{S}_s$, where $\mathbb{S}_s\subset \R^{p\tau}$ is the set of block $s$-sparse vectors.

We remark that equation \eqref{cs} has a unique solution in $\R^n\times \mathbb{S}_s$ if and only if the CPS defined in \eqref{system} is $2s$-sparse observable, that is, by removing any subset of $2s$ sensors, the attack-free system remains observable; see, e.g., \cite[Theorem 3.2]{sho16}. In other terms, $2s$-sparse observability is equivalent to the injectivity of the map $f:\R^n\times \mathbb{S}_s\mapsto \R^{p\tau}$ defined as $f(x,a)=\Obs x +a$. This condition is necessary condition for secure state estimation, independently from the estimation algorithm.

If on the one hand the prior information on the block-sparsity pattern of $\atrue$ can improve the state estimation, on the other hand it ties the solution to constant attack supports and it originates more complex recovery algorithms.
For these motivations, in this work  we propose a sparse optimization approach that neglects the possible block-sparsity.

 %Moreover, as explained in \cite[Sec. III.A]{sho16},  sparse observability is equivalent to the concept of strong observability of DT LTI systems with unknown input, see, e.g., \cite{sun11}.
%In particular, the $2s$-sparse observability implies $p>2s$.

Beyond the seminal contributions briefly described in this section, we remark that substantial recent work addresses specific aspects of SSE of CPSs of the kind \eqref{system}. For example, \cite{mao22} and \cite{lu23} develop decomposition techniques, \cite{zha23} propose a data-driven approach; \cite{lei23,mao23} focus on distributed algorithms for SSE.
%
% in thi
% % Although this paper focsus on general and foundational models for CPSs under sensor attacks, in order to better motivate our study, we...
% %
% % , respectovely based on the generalized eigenspaces of the system matrix $A$ and by adding information about the measurement matrix $C$. Decomposition techniques consist in state transformations that allow to identify portions of the state than one can recover in polynomial time. \cite{zha23} propose a data-driven approach is proposed for the detection/identification of sparse sensor attacks in linear CPSs with unkown parameters. \cite{lu23b} introduce the use of side information either from other transmitters or from the sensors themselves.
% %
% % Finally, several recent works introduce distributed algorithms for SSR, see, e.g., \cite{lei23,mao23,shi23}.
% %

\section{Lasso approach}\label{sec:PA}
The proposed Lasso formulation for problem \eqref{cs} is
\begin{equation}\label{lasso}
(\xstar,\astar)=\argmin{x\in\R^n,a\in\R^{p\tau}} \frac{1}{2}\left\|y - \Obs x -a\right\|_2^2+\lambda \|a\|_1
\end{equation}
where $\lambda>0$. Problem \eqref{lasso} is a partial Lasso because only a part of the vector to estimate is sparse, i.e., we apply the $\ell_1$ regularization only on $a$.

We notice that constrained counterparts of Lasso are the bases for the $\ell_1$-based, block-sparsity decoders proposed by \cite{faw14}, for the noise-free case, and by  \cite{paj17}, in the presence of bounded measurement noise. The  Lasso formulation in \eqref{lasso} may envisage the presence of measurement noise as well; moreover, we can resort to several effective algorithms for unconstrained optimization to solve it.

In this work, we consider the iterative soft thresholding algorithm (ISTA) proposed by \cite{dau04}, which is a proximal gradient algorithm. ISTA iteration consists of a gradient step and a componentwise soft thresholding operation, defined by $\soft_{\nu\lambda}[w]=w-\nu\lambda\sign(w) $ if $|w|\geq \nu\lambda$, and 0 otherwise; $\nu>0$ is the gradient step size. Accelerated versions of ISTA are available, see, e.g.,  \cite{bec09} and \cite{fox23}. Its simple structure allows us to build a state observer upon it, as illustrated in Sec. \ref{sec:OB}.
%ISTA is similar to ETPG, but ETPG performs an event-triggered (not at each iteration) hard thresholding, which saves the $s$ largest components. In particular, ETPG requires to know the sparsity level and time-invariant sparsity pattern.
%
%Our interest in ISTA is motivated also by the possibility of reformulating it as a sparse observer
%
%To accelerate ISTA, we can implement its fast version FISTA proposed by \cite{bec09} or non-convex methods proposed by \cite{fox23}.
%
\subsection{Analysis of the irrepresentable condition}
An interesting feature of Lasso is that there is a tight condition, denoted as ``irrepresentable'', that guarantees the recovery of the correct support in the noise-free case, see \cite{fuc04}. Extensions to the noisy measurements are possible as well, see, e.g.,   \cite{fuc05}.

In a nutshell, given a classic Lasso $\frac{1}{2}\left\|y - Qz\right\|_2^2+\lambda \|z\|_1$, the irrepresentable condition states that the columns of the sensing matrix $Q$ on the true support must be ``sufficiently orthogonal'' to the columns outside the support.
The irrepresentable condition cannot be priorly checked, because it involves the knowledge of the support; however, it provides interesting insights on the necessary features of $Q$ to recover the support through Lasso; see, e.g., \cite{zha06,has15book,fox20irr}.

In the considered SSE problem, the sensing matrix $\begin{pmatrix}
                                        \Obs&I
                                       \end{pmatrix}$ has a peculiar structure, with an identity matrix in its right part.
In this section, we perform an irrepresentable condition analysis that takes into account this structure. We focus on the attack support recovery because its correct estimation allows us to recover the state from the safe sensors.

Let $\su$ be the support of $\atrue$, and $|\su|=h \ll p\tau$. If the attack support is time-invariant, $h=s\tau$.
In the rest of the paper, we use the following notation: $\Is\in\{0,1\}^{p\tau,h}$ is the submatrix of $I$ with columns indexed in $\su$. $\Obs_{\su}\in\R^{h,n}$ and $\Obs_{\suc}\in\R^{p\tau-h,n}$ are the submatrices of $\Obs$ with \emph{rows} in $\su$ and in $\suc$, respectively.
$\Obs_{\suc}^{\dagger\top}=\Obs_{\suc}[\Obs_{\suc}^\top\Obs_{\suc}]^{-1}\in\R^{p\tau-h,n}$
is the right pseudo-inverse of $\Obs_{\suc}^\top$. Finally, we denote by $\|\cdot\|_{\infty}$ the $\ell_{\infty}$ matrix norm, defined as $\|M\|_{\infty}=\max_{i}\sum_{j} |M_{i,j}|$ for any matrix $M$. If $M$ is a row vector, the $\ell_{\infty}$ matrix norm corresponds to the $\ell_1$ vector norm.

The following result holds.
\begin{theorem}
Let us assume that $\begin{pmatrix}
                     \Obs& \Is
                    \end{pmatrix}\in\R^{p\tau,n+h}
$ is full rank. Lasso is successful, i.e., by solving it we identify the attack support, if and only if
\begin{equation}\label{irr_stretta}
 \left\|\Obs_{\suc}^{\dagger\top}\Obs_{\su}^\top\sign(\atrue_{\su})\right\|_{\infty}<1
\end{equation}
provided that $\lambda>0$ is sufficiently small.
As a consequence, Lasso is successful if
\begin{equation}\label{irr_suff}
 \left\|\Obs_{\suc}^{\dagger\top}\Obs_{\su}^\top\right\|_{\infty}<1.
\end{equation}

\end{theorem}
%We assess the value of $\lambda$ in the proof.
\begin{proof}
$(\xstar,\astar)$ is a global minimum of \eqref{lasso} if and only if it fulfills the zero-subgradient condition: there exists $\zeta \in \partial \|\astar\|_1$ such that
\begin{equation}\label{subg}
\begin{pmatrix}
\Obs^\top \\ I
\end{pmatrix}\left(\Obs \xstar +\astar -y \right)+ \lambda\begin{pmatrix}
                                       0\\ \zeta
                                      \end{pmatrix}=0
\end{equation}
where $\partial \|\astar\|_1\subset \R^{p\tau}$ is the subdifferential of $\|\astar\|_1$: for each $i\in\su$, $\zeta_i=\sign(\astar_i)$, while $\|\zeta_{\suc}\|_{\infty}\leq 1$; see \cite{fuc04} for details. In particular, if $\|\zeta_{\suc}\|_{\infty}< 1$, then $(\xstar,\astar)$ is a strict minimum.

We construct the candidate $(\xstar,\astar) \in\R^{n+p\tau}$ as follows: we set $\astar_{\suc}=0$, while
 \begin{equation}
 (\xstar,\astar_{\su})= \argmin{x\in\R^n,a_{\su}\in\R^{h}} \frac{1}{2}\left\|y - \Obs x -I_{\su}a_{\su}\right\|_2^2+\lambda \|a_{\su}\|_1
 \end{equation}
 which is a Lasso restricted on the non-zero components of the true vector.
 %The explicit solution is $(\xstar,\astar_{\su})=$
%
 The so-built $(\xstar,\astar)$ is a solution of Lasso \eqref{lasso} if it satisfies \eqref{subg}; let us verify under which conditions this occurs.

 By distinguishing the zero-subgradient equations on $(\xstar,\astar_{\su})$ and on $\astar_{\suc}$ and by recalling that $y=\Obs\xtrue + I_{\su}\atrue_{\su}$, we have
\begin{equation}\label{onSx}
\begin{pmatrix}
\Obs& I_{\su}
\end{pmatrix}^{\top}\begin{pmatrix}\Obs & \Is \end{pmatrix}\begin{pmatrix}\xstar-\xtrue \\ \astar_{\su}-\atrue_{\su} \end{pmatrix}+ \lambda\begin{pmatrix}0\\
 \sign(\astar_{\su})\end{pmatrix}=0
% \begin{pmatrix}
% \Obs^{\top} \\ I_{\su}^{\top}
% \end{pmatrix}\left(\Obs x +I_{\su} \as -y \right)+ \lambda\begin{pmatrix}0\\
% \sign(\as)\end{pmatrix}=0
\end{equation}
\begin{equation}\label{onSc} I_{\suc}^{\top}\begin{pmatrix}\Obs & \Is \end{pmatrix}\begin{pmatrix}\xstar-\xtrue \\ \astar_{\su}-\atrue_{\su} \end{pmatrix}+ \lambda
\zeta_{\suc}=0.
\end{equation}
%where $\|\zeta_{\suc}\|_{\infty}\leq 1$.
%
%
We notice that
\begin{equation}\label{block} \begin{pmatrix}
\Obs & I_{\su}
\end{pmatrix}^{\top}\begin{pmatrix}\Obs & \Is \end{pmatrix}=\begin{pmatrix}\Obs^\top\Obs& \Obs_{\su}^\top\\ \Obs_{\su}& I_h \end{pmatrix}\end{equation} where  $I_h\in\{0,1\}^{h,h}$ is the identity matrix of dimension $h$.
Then, from \eqref{onSx} and \eqref{block}, we compute
\begin{equation}\label{onS}
 \begin{pmatrix}\xstar-\xtrue \\ \astar_{\su}-\atrue_{\su} \end{pmatrix}=-\lambda \begin{pmatrix}\Obs^\top\Obs& \Obs_{\su}^\top\\ \Obs_{\su}& I_h \end{pmatrix}^{-1}\begin{pmatrix}0\\
\sign(\astar_{\su})\end{pmatrix}.
\end{equation}
The inverse matrix exists because the rows of $\begin{pmatrix}\Obs & \Is \end{pmatrix}^\top\in\R^{n+h,p\tau}$, $n+h<p\tau$, are linearly independent by assumption.
Moreover, from \eqref{onS}, the difference $\astar_{\su}-\atrue_{\su}$ is proportional to $\lambda$; therefore, if $\lambda$ is sufficiently small, $\sign(\astar_{\su})=\sign(\atrue_{\su})$. We remark that, in the noise-free case, we can design $\lambda$ as  arbitrarily small without loss of generality.

By replacing \eqref{onS} in \eqref{onSc}, we obtain;
\begin{equation}
%\begin{split}
  \zeta_{\suc}=I_{\suc}^{\top}\begin{pmatrix}\Obs & \Is \end{pmatrix} \begin{pmatrix}\Obs^\top\Obs& \Obs_{\su}^\top\\ \Obs_{\su}& I_h \end{pmatrix}^{-1}\begin{pmatrix}0\\
\sign(\atrue_{\su})\end{pmatrix}.\\
%&=I_{\su}^{\top}\begin{pmatrix}\Obs & \Is \end{pmatrix} \left[\begin{pmatrix}\Obs^\top \\ \Is^\top \end{pmatrix}\begin{pmatrix}\Obs & \Is \end{pmatrix}\right]^{-1}\begin{pmatrix}0\\
%\sign(\as)\end{pmatrix}.
%\end{split}
\end{equation}
Now, $(\xstar,\astar)$ is the unique minimum of Lasso if $\|\zeta_{\suc}\|_{\infty}< 1$. Therefore, we study the inequality
\begin{equation}\label{conditiontoverify}
 \left\|I_{\suc}^{\top}\begin{pmatrix}\Obs & \Is \end{pmatrix} \begin{pmatrix}\Obs^\top\Obs& \Obs_{\su}^\top\\ \Obs_{\su}& I_h \end{pmatrix}^{-1}\begin{pmatrix}0\\
\sign(\atrue_{\su})\end{pmatrix}\right\|_{\infty}<1.
\end{equation}
Since $\Isc^\top\Is=0$,
\begin{equation}\label{ort}
 I_{\suc}^{\top}\begin{pmatrix}\Obs & \Is \end{pmatrix}= \begin{pmatrix}\Obs_{\suc} & 0 \end{pmatrix}.
\end{equation}
%where $\Obs_{\suc}\in\R^{p\tau-h,n}$ contains the \emph{rows} of $\Obs$ in $\suc$.
%
From Schur's complement arguments,
\begin{equation}\label{inv}
\begin{pmatrix}\Obs^\top\Obs& \Obs_{\su}^\top\\ \Obs_{\su}& I_h \end{pmatrix}^{-1}=\begin{pmatrix}\Omega_1& \Omega_2\\ \Omega_3& \Omega_4 \end{pmatrix}
\end{equation}
where
\begin{equation}\label{o2}
 \Omega_2=-[\Obs^\top\Obs-\Obs_{\su}^\top\Obs_{\su}]^{-1}\Obs_{\su}=-[\Obs_{\suc}^\top\Obs_{\suc}]^{-1}\Obs_{\su}^\top.
\end{equation}
By applying \eqref{ort},\eqref{inv} and \eqref{o2} to \eqref{conditiontoverify},
\begin{equation}
 \begin{pmatrix}\Obs_{\suc} & 0 \end{pmatrix} \begin{pmatrix}\Omega_1& \Omega_2\\ \Omega_3& \Omega_4 \end{pmatrix}\begin{pmatrix}0\\
\sign(\atrue_{\su})\end{pmatrix}=-\Obs_{\suc}^{\dagger\top}\Obs_{\su}^\top\sign(\atrue_{\su}).
\end{equation}
%where
%\begin{equation}
%\Obs_{\suc}^{\dagger\top}=\Obs_{\suc}[\Obs_{\suc}^\top\Obs_{\suc}]^{-1}\in\R^{p\tau-h,n}
%\end{equation}
%is the right pseudo-inverse of $\Obs_{\suc}^\top$.
In conclusion, condition \eqref{conditiontoverify} is equivalent to \eqref{irr_stretta}.

%\begin{equation}
 %$$\left\|\Obs_{\suc}^{\dagger\top}\Obs_{\su}^\top\sign(\atrue_{\su})\right\|_{\infty}<1.$
%\end{equation}
Since
 %\begin{equation}
 %begin{split}
   $\left\|\Obs_{\suc}^{\dagger\top}\Obs_{\su}^\top\sign(\atrue_{\su})\right\|_{\infty}\leq \left\|\Obs_{\suc}^{\dagger\top}\Obs_{\su}^\top\right\|_{\infty}$
%   %\leq \left\|\Obs_{\suc}^{\dagger\top}\Obs_{\su}^\top\right\|_{\infty}\left\|\sign(\atrue_{\su})\right\|_{\infty}\\
%  % &~~~~~~~~~~~~~\leq \left\|\Obs_{\suc}^{\dagger\top}\Obs_{\su}^\top\right\|_{\infty}=\|\Obs_{\su}\Obs^\dagger_{\suc}\|_1
% \end{split}\end{equation}
 then \eqref{irr_suff} is sufficient for a successful Lasso.
%\QED
\end{proof}

% [SULL'INTERPRETAZIONE MI PERDO]
%As explained in \cite[Sec. 11.4.1]{has15book}, in a generic Lasso, we can interpreted the irrepresentable condition as follows: the irrepresentable condition holds if the columns of the sensing matrix in the support and the ones outside the support are sufficiently ``orthogonal'' or incoherent.
%A difference between \eqref{irr_stretta} and \eqref{irr_suff} is that \eqref{irr_stretta} depends on the signs of the attacks; instead, \eqref{irr_suff} only depends on the observability matrix.

Since $\irr:=\left\|\Obs_{\suc}^{\dagger\top}\Obs_{\su}^\top\right\|_{\infty}=\max_{j}\left\|\Obs_{\su}\left(\Obs^\dagger_{\suc}\right)_j\right\|_1$, a qualitative interpretation of \eqref{irr_suff} is as follows: the term $\irr$  must be small, i.e., the rows of $\Obs_{\su}$ must be ``sufficiently orthogonal'' to the columns of $\Obs_{\suc}^{\dagger}$.
This implies some observations. Since $\xtrue=\Obs_{\suc}^{\dagger} y_{\suc}$, then
\begin{equation}
\Obs_{\su}\xtrue=\Obs_{\su}\Obs_{\suc}^{\dagger} y_{\suc}.
\end{equation}
In particular, $\Obs_{\su}\xtrue=0$ would be the ideal case to identify the attacks, which would be directly observable from $y_{\su}=a_{\su}$. Nevertheless, this is not realistic, in particular because it depends on the specific initial state $\xtrue$. However,
\begin{equation}
% \begin{split}
% \left\|\Obs_{\su}\xtrue\right\|_1&=\left\|\sum_{j=1,\dots,(p-h)\tau}\Obs_{\su}\left(\Obs_{\suc}^{\dagger}\right)_j \left(y_{\suc}\right)_j\right\|_1\\
% &\leq\sum_{j=1,\dots,(p-h)\tau}\left\|\Obs_{\su}\left(\Obs_{\suc}^{\dagger}\right)_j\right\|_1 | \left(y_{\suc}\right)_j|\\
% &\leq \irr\left\| y_{\suc}\right\|_1
%\end{split}
 \left\|\Obs_{\su}\xtrue\right\|_1\leq \rho\left\| y_{\suc}\right\|_1
\end{equation}
that is, $\irr$ controls the ``orthogonality'' between the rows of $\Obs$ indexed in $\su$ and $\xtrue$: a small  irrepresentable term $\irr$ implies a small energy of $\Obs_{\su}\xtrue$, which, in turn, implies that the attacks are more  exposed to identification.
%
%
%\subsection{Illustrative examples}

Differently from previous conditions considered in the literature, see for example Theorem 4.4 in \cite{sho16}, the irrepresentable condition well captures the ability of CPS to identify attacks. We illustrate this point with a simple illustrative example.
\begin{example}
Let us consider a simple static model with $A=I$, $n=1$ and $\Obs=C=\alpha\begin{pmatrix}
                                                                  1&1&1
                                                                 \end{pmatrix}^\top$ for any $\alpha\in\R$, i.e., we have 3 equivalent sensors. We assume that one of them is under attack. If we know that $h=1$, the identification of the attack is trivial: the sensor that provides a different measurement is clearly under attack and we can recover the state from the other two sensors.
                                                                 The irrepresentable condition well captures this resilience; in fact, we have $\Obs_{\suc}^{\dagger\top}=\frac{1}{2\alpha}\begin{pmatrix}
                                                                  1&1
                                                                 \end{pmatrix}^\top$ and $\Obs_{\su}^\top=\alpha$, thus $\left\|\Obs_{\suc}^{\dagger\top}\Obs_{\su}^\top\right\|_{\infty}=\frac{1}{2}<1$ for any $\alpha\in\R$ and Lasso is successful.

In contrast, we notice that condition 2 of Theorem 4.4 in \cite{sho16} that guarantees the convergence of ETPG is not satisfied. In fact, the maximum eigenvalue of $\begin{pmatrix}
                                               \Obs&I
                                              \end{pmatrix}^\top\begin{pmatrix}
                                               \Obs&I
                                              \end{pmatrix}
$ is $q=3\alpha^2+1$, while the minimum eigenvalue on all the possible subsets of 2 sensors is $r=\frac{3\alpha^2+1-\sqrt{9\alpha^4+2\alpha^2+1}}{2}$. Then, $r<\frac{4}{9}q$, which contradicts condition 2. %In conclusion, differently from the irrepresentable condition, the theory by \cite{sho16}  does not capture the resilience of this CPS.
\end{example}
%
% \begin{example}
% We consider the problem of controlling an unmanned ground vehicle (UGV) studied in \cite[Sec. V.B]{sho17}. This is a system with two states, namely position and velocity of the vehicle. A safe GPS sensor measures the position, while two vulnerable motor encoders measure the velocity. We consider the parameter values and the discretization setting as proposed by \cite{sho17}, which are publicly valuable at \cite{code}.  In turn, each motor encoder may be subject to attacks, while the GPS sensor is attack-free. In this setting, the system is shown to be 2-observable. On the other hand, the condition 2 of \cite[Theorem 4.4)]{sho16} does not hold, therefore ETPG is not guaranteed to converge.
%
% In this example, the irrepresentable condition \eqref{irr_suff} does not hold, since $\irr = 1.006$; then, Lasso is not guaranteed to be successful, while condition \eqref{irr_stretta} holds with value 0.006 if $\sign(\atrue_{\su})=\pm\begin{pmatrix}
%                                                                                1&-1
%                                                                               \end{pmatrix}^\top
% $, i.e., if at two successive time steps, the attack has switched sign. From numerical simulations, we double-check that Lasso correctly locate the attack.
% \end{example}

\section{Sparse soft observer for online SSE}\label{sec:OB}
In this section, we move towards recursive, online SSE. We consider Problem \ref{p1} in a dynamic perspective: we aim at estimating the current state, or a delayed version, using the last $p\tau$ measurements, and, at each $k$ we include the new measurements and we discard the oldest  ones. If $\tau=1$, this an online (not delayed) SSE. This calls for fast recursive online algorithms.

\cite{sho16} address this problem by developing a recursive version of ETPG, named ETPL. At time $k$, instead of running ETPG to convergence to estimate $x(k-\tau+1)$, ETPL runs some steps of the algorithm, then it moves to step $k+1$ and  it takes new measurements, in the philosophy of a Luenberger observer.
As an alternative, we develop an online version of ISTA, that we name sparse soft observer and that we summarize in Alg. \ref{alg}. We use the following notation: $\mathbf{a}(k)=(a(k-\tau+1)^\top,\dots,a(k)^\top)^\top$,
 $\mathbf{y}(k)=(y(k-\tau+1)^\top,\dots,y(k)^\top)^\top$.

\begin{algorithm}
  \renewcommand{\algorithmicrequire}{\textbf{Input:}}
    \renewcommand{\algorithmicensure}{\textbf{Output:}}
\setstretch{1.2}
 \caption{Sparse soft observer}
  \begin{algorithmic}[1]
  \REQUIRE $\tau\leq n$, $\lambda>0$, $\nu>0$, $A$, $\Obs$, $\mathbf{y}(k)$
 \ENSURE $\begin{pmatrix}
    \xhat(k)\\
    \ahat(k)
   \end{pmatrix} = $ estimate of $\begin{pmatrix}
    x(k-\tau+1)\\
    \mathbf{a}(k)
   \end{pmatrix}$
%Given $y(k)=Cx(k)+a$ and  $(\xhat(k),\ahat(k))^{\top}$ computed before time $k$
 \FORALL{$k=\tau-1,\tau,\dots$}
\STATE Measurements and estimated measurements update
   \begin{equation}
   \begin{split}
    \mathbf{y}(k)&=\Obs x(k-\tau+1)+\mathbf{a}(k)\\
  \yhat(k)&=\Obs\xhat(k)+\ahat(k)\\
    \end{split}
   \end{equation}
\STATE ISTA step: gradient step + soft thresholding
 \begin{equation}
   \begin{split}
 \begin{pmatrix}
     \xhat^+\\
     \ahat^+
    \end{pmatrix}
= &\begin{pmatrix}
    \xhat(k)\\
    \ahat(k)
   \end{pmatrix}- \nu \begin{pmatrix}
                                               \Obs&I
                                              \end{pmatrix}^{\top}[\yhat(k)-\mathbf{y}(k)]
                                              \end{split}
                                              \end{equation}
\begin{equation}
\begin{split}
%\xhat(k+1)&=A\xhat^+\\
\ahat(k+1)&=\soft_{\nu\lambda}\left[\ahat^+\right]
   \end{split}
   \end{equation}
% = &\begin{pmatrix}
%     \xhat(k)\\
%     \ahat(k)
%    \end{pmatrix}- \nu \begin{pmatrix}
%                                                \Obs&I
%                                               \end{pmatrix}^{\top}\begin{pmatrix}
%                                                \Obs&I
%                                               \end{pmatrix}\left(\begin{pmatrix}
%     \xhat(k)\\
%     \ahat(k)
%    \end{pmatrix}-\begin{pmatrix}
%     x(k-\tau+1)\\
%     \mathbf{a}(k)
%    \end{pmatrix}\right)\\
% %   &= \begin{pmatrix}
% %     \xhat(k)\\
% %     \ahat(k)
% %    \end{pmatrix}- \tau \begin{pmatrix}
% %                                                \Obs&I
% %                                               \end{pmatrix}^{\top}\begin{pmatrix}
% %                                                \Obs&I
% %                                               \end{pmatrix}[\begin{pmatrix}
% %     x(k-\tau+1)\\
% %     \mathbf{a}(k)
% %    \end{pmatrix}-\begin{pmatrix}
% %     x(k)\\
% %     \ahat(k)
% %    \end{pmatrix}]\\
\STATE State update
\begin{equation}
\begin{split}
\xhat(k+1)&=A\xhat^+\\
%\ahat(k+1)&=\soft_{\nu\lambda}\left[\ahat^+\right]
   \end{split}
   \end{equation}
\ENDFOR
 \end{algorithmic}\label{alg}
\end{algorithm}

The basic idea is to run one ISTA step, as described in Sec. \ref{sec:PA},  at each time instant $k$. Then, we update the estimate of the state by leveraging the knowledge of $A$, as in Luenberger observer. Step 3. in Alg. \ref{alg} can be repeated more than one time to enhance the estimation, but we do not run the complete algorithm to converge to the Lasso solution. Moreover, in case of online estimation, i.e., $\tau=1$, the current $p$ measurements are expected to be insufficient to have a successful Lasso; therefore, the observer approach is necessary.

\section{Numerical results}\label{sec:NR}
In this section, we propose some numerical results to investigate the performance of the proposed Lasso approach and sparse soft observer, in terms of state estimation accuracy and execution time. We perform all the simulations in MATLAB R2023 on a processor i7 @ 1.80 GHz $\times$ 8, with 16 GB of RAM.
%Il ciccio: i7-10875H CPU @ 2.30GHz × 16, with 32 GB of RAM
\subsection{Lasso approach}
We test Lasso, solved by FISTA, see \cite{bec09}, on random, synthetic CPSs and we compare it to ETPG by \cite{sho16} and Imhotep-SMT by \cite{sho17}; the code for this last one is taken at \cite{code}. We generate the elements $A$ and $C$ independently, according to a standard normal distribution; then we normalize $A$ to guarantee stability.
We assume that the attack support is time-invariant and cardinality $s$ and we generate it uniformly at random. The initial state $x(0)$ has uniformly distributed components with magnitude in $[2,3]$. The attacks have magnitude in $[4,5]$, which is sufficiently large to sabotage the state estimation, but not enough large to produce clear, plainly detectable outliers in the measurements.

\begin{figure}
\centering
\includegraphics[width=0.45\textwidth]{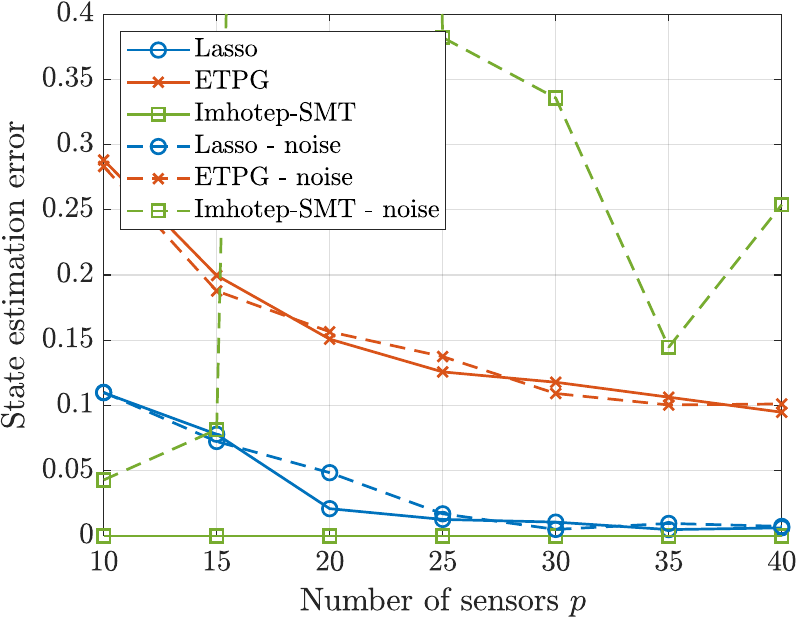}\hspace{0.8cm}
\includegraphics[width=0.45\textwidth]{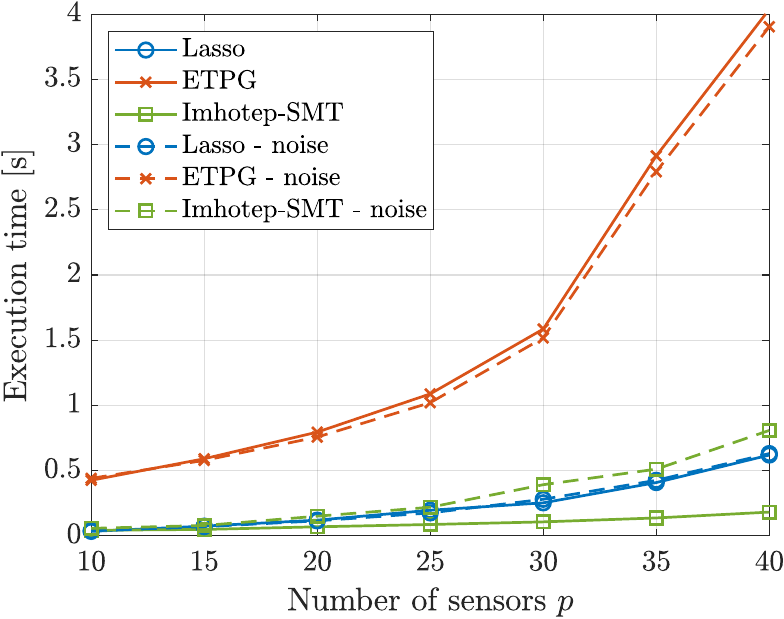}
\caption{Lasso vs ETPG vs Imhotep-SMT, $n=20$, $s=p/5$, noise-free and with noise bound $10^{-4}$}\label{fig:p}
\end{figure}

\begin{figure}
\centering
\includegraphics[width=0.46\textwidth]{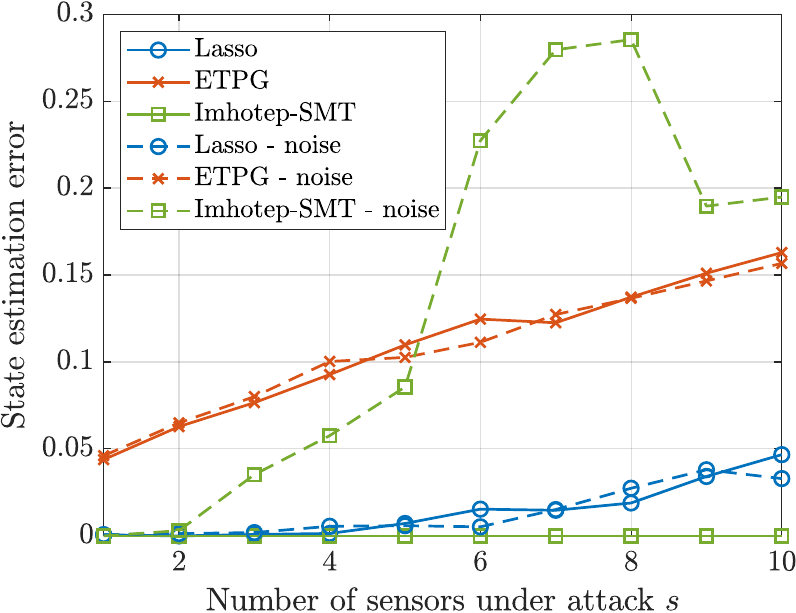}\hspace{0.8cm}
\includegraphics[width=0.46\textwidth]{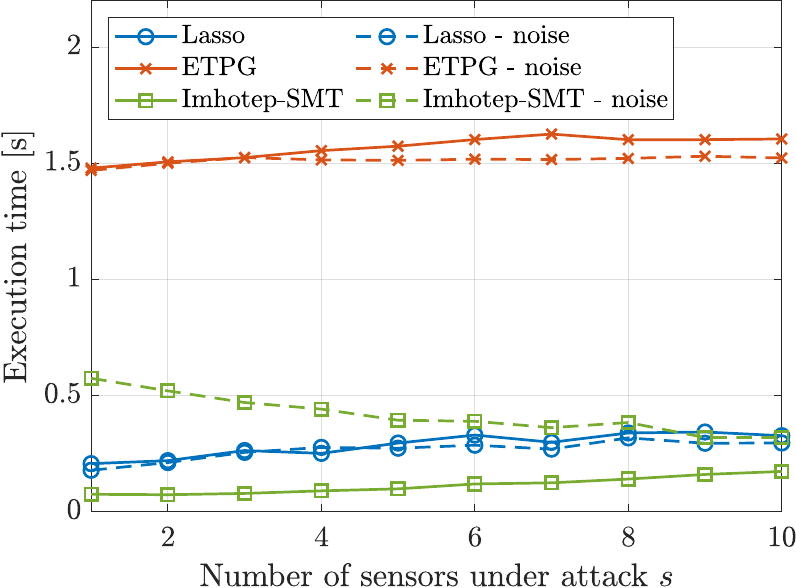}
\caption{Lasso vs ETPG vs Imhotep-SMT, $n=20$, $p=30$, noise-free and with noise bound $10^{-4}$}\label{fig:s}
\end{figure}

In the first experiment, we vary $p$, while $n=20$ and $s=\frac{p}{5}$; in the second experiment, we vary $s$, while $n=20$ and $p=30$.
%We perform two classes of experiments: in each of them, we vary one parameter among $n$, $p$ and $s$ and we keep constant the others.
We perform 50 runs for each experiment; we depict the averages results in Fig. \ref{fig:p} and in Fig. \ref{fig:s}.
We assess the accuracy in terms of state estimation error $\|\xhat-\xtrue\|_2/\|\xtrue\|_2$.  We consider either noise-free and noisy measurements, i.e., $y(k)=Cx(k)+a(k)+\eta(k)$, where $\eta(k)\in\R^p$ is a uniformly random, bounded noise with $\|\eta(k)\|_{\infty}\leq 10^{-4}$.

As we can see in Fig. \ref{fig:p} and in Fig. \ref{fig:s}, in this experiment, Lasso outperforms ETPG both in accuracy and run time. Since we consider small/medium dimensions, Imhotep-SMT is the best approach to provide the exact solution in fast time, in the noise-free case; nevertheless, its accuracy is critical in the presence of small measurement noise. In contrast, Lasso and ETPG are robust to small noise, with a very slight degradation in the estimation accuracy. %In conclusion, Lasso is preferable with respect to ETPG and Imhotep-SMT in small/medium dimension CPSs, in the presence of measuremement noise. We leave to future work large dimensional CPSs, where the combinatorial nature of Imhotep-SMT may become computationally critical.

% \begin{figure*}\label{fig:n}
% \centering
% \includegraphics[width=0.4\textwidth]{compact_n_est}\hspace{0.8cm}
% \includegraphics[width=0.4\textwidth]{compact_n_time}
% \caption{Lasso vs ETPG vs Imhotep-SMT, $p=30$, $s=5$, noise-free and with noise bound $10^{-4}$}
% \end{figure*}
\subsection{Sparse soft observer}
In this section, we test the proposed sparse soft observer for recursive and online SSR. We perform 100  random runs and we show the average results, in terms of time evolution of the state estimation error $\|\xhat-\xtrue\|_2/\|\xtrue\|_2$ and support error, defined as $\sum_j |\mathbf{1}(\ahat_j \neq 0)- \mathbf{1}(\atrue_j \neq 0)|$, where $\mathbf{1}(v)=1$ if $v$ is true and 0 otherwise.

We generate noise-free dynamical models as in the previous experiments, with $n=10$, $p=15$, $s=3$, and we run each system for 300 time steps. We implement the proposed sparse soft observer and, for comparison, ETPL by \cite{sho16}.
\begin{figure}
\centering
\includegraphics[width=0.46\textwidth]{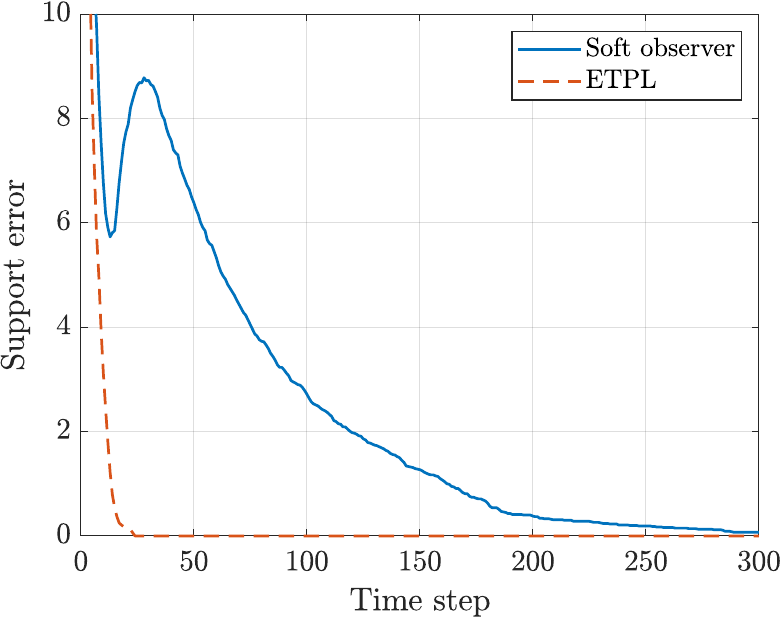}
\includegraphics[width=0.47\textwidth]{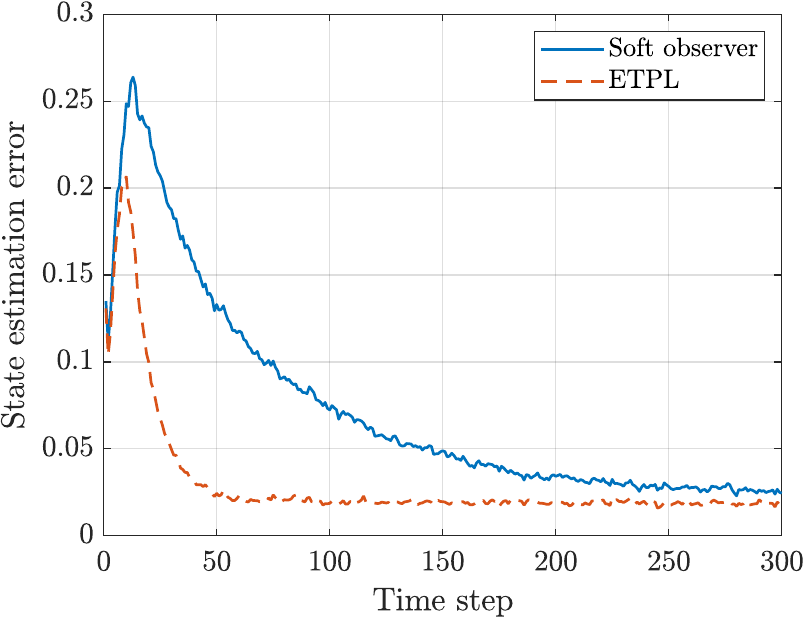}
\caption{Sparse soft observer vs ETPL; $n=10$, $p=15$, $s=3$, $\tau=n$}\label{fig:track1}
\end{figure}
\begin{figure}
\centering
\includegraphics[width=0.46\textwidth]{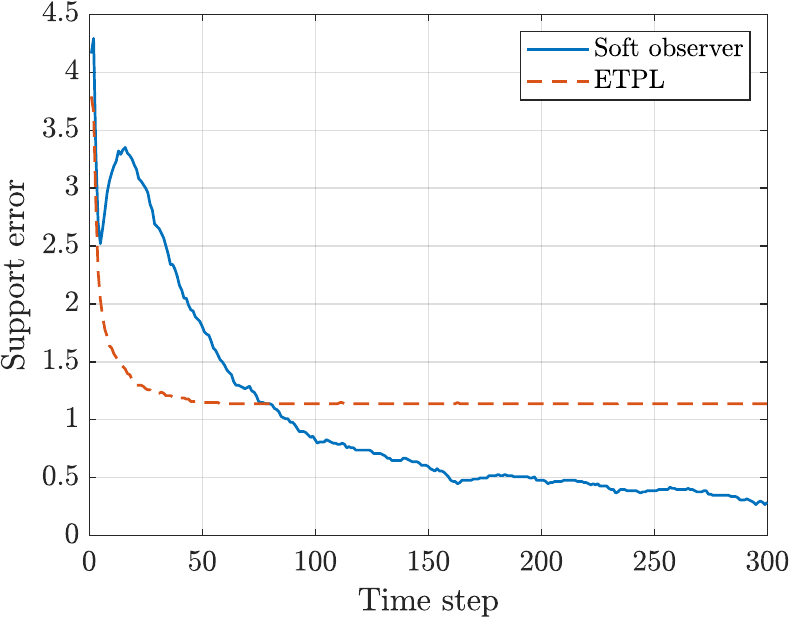}
\includegraphics[width=0.46\textwidth]{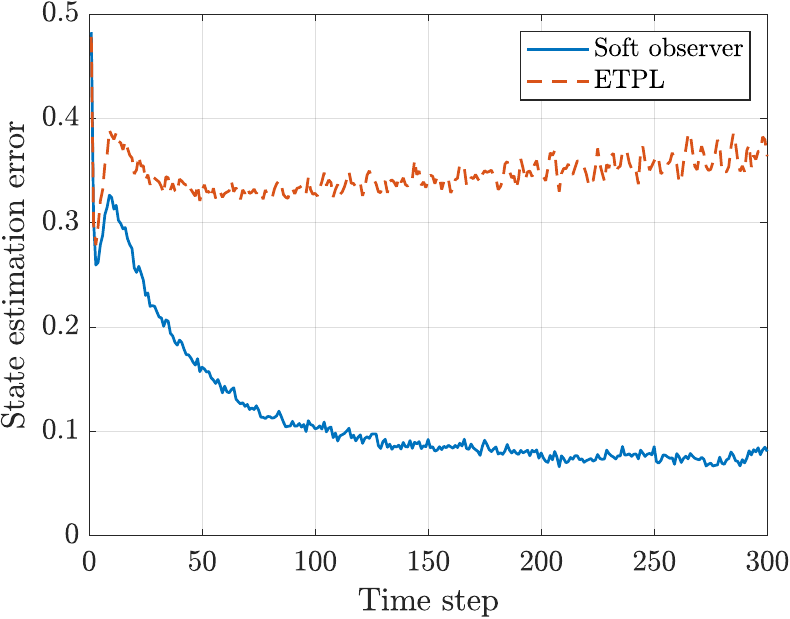}
\caption{Sparse soft observer vs ETPL; $n=10$, $p=15$, $s=3$, $\tau=1$}\label{fig:track2}
\end{figure}
We show the results in Fig. \ref{fig:track1} and in Fig. \ref{fig:track2}. In Fig. \ref{fig:track1}, we set $\tau=n$, that is, at each time step $k$ we use the previous $p\tau$ measurements to estimate the delayed $x(k-\tau+1)$. Then, at each $k$, we remove the oldest $p$ measurements to introduce the new $p$ ones. This is the case considered also by \cite{sho16}. Instead, in Fig. \ref{fig:track2}, we consider $\tau=1$, that is, the algorithms use the current set of measurements $y(k)$ to provide an online estimate of $x(k)$. In the sparse soft observer, we repeat the step 3 of Alg. \ref{alg}. for $5\tau$ times.
We can see that ETPL is more accurate for $\tau=n$, in particular it converges more quickly than the sparse soft observer to the correct attack support. However, the execution time at each $k$ is $2\cdot 10^{-3}$ seconds for ETPL and $4\cdot 10^{-4}$ seconds for the sparse soft observer, so the latter is adaptable to fast time scales.
On the other hand, for $\tau=1$, the sparse observer is more accurate and ETPL does not always converge to the right support. In this case, the execution times are $7\cdot 10^{-5}$ seconds for ETPL and $4\cdot 10^{-6}$ seconds for the sparse soft observer.
\section{Conclusions}\label{sec:C}
We propose a Lasso approach for secure state estimation in cyber-physical systems under sparse sensor attacks. We analyse the properties of Lasso to identify  the attack and, as a consequence, to recover the state. Furthermore, by starting from the iterative soft thresholding algorithm for Lasso, we develop a sparse soft observer to perform online estimation. Through numerical results, we show that the proposed Lasso approach is valuable with respect to state-of-the-art methods, although it exploits less information, e.g., on the sparsity pattern. Moreover, in our experiments, the sparse soft observer converges to sufficiently accurate solutions, with a reduced execution time.
Future work includes the extension of the analysis to noisy models and the study of the convergence of the sparse soft observer.

%AD-ISTA, a variant of ISTA developed by applying the proximal gradient method to Log-Lasso. AD-ISTA converges in less iterations with respect to ISTA, FISTA and ADMM, thanks to an adaptive shrinkage hyperparameter, that limits the increase of the $\ell_1$-norm during the first phase. Moreover, by applying  the principles of FISTA, we also propose the accelerated version AD-FISTA. Through numerical experiments, we verify that AD-ISTA is faster than the state-of-the-art algorithms for Lasso and that we obtain a further acceleration with AD-FISTA. Possible extensions of this work include the rigorous proof of the convergence rate and the generalization to sparse optimization problems different from Lasso.
%
%\begin{ack}
%The authors would like to thank Prof. Jean-Bernard Lasserre for his constructive comments on the manuscript.
%\end{ack}
%
\bibliographystyle{plain}
\bibliography{CPS_refs_nov2023}             % bib file to produce the bibliography
\end{document}